\documentstyle{article}
\begin{document}
\newtheorem{thm}{Theorem}[section]
\newtheorem{cor}{Corollary}[section]
\newtheorem{prop}{Proposition}[section]

\begin{center}
{\Large A representation formula for non-conformal 

\vspace{2mm}

harmonic surfaces in $R^3$} 

\vspace{5mm}

Bart\ Dioos and Makoto\ Sakaki 
\end{center}

\vspace{5mm}

Abstract. We discuss non-conformal harmonic surfaces in $R^3$ with prescribed ($\pm$)transforms, and we get a representation formula for non-conformal harmonic surfaces in $R^3$. 

\vspace{5mm}



\section{Introduction}

As a national generalization of minimal surfaces, harmonic surfaces have been studied (see \cite{AL}, \cite{C}, \cite{DVV}, \cite{K}, \cite{KM}, \cite{S2} and references therein). In this paper, we study non-conformal harmonic surfaces in the Euclidean $3$-space $R^3$, from the view point of transforms and representation. 

Classically, it is well known that the Gauss map of a minimal surface in $R^3$ gives a holomorphic map to the $2$-sphere $S^2$ (cf. \cite{O}). For a minimal surface in the $3$-sphere $S^3$, the Gauss map gives a minimal surface in $S^3$ possibly with singularities (cf. \cite{L}). This is the polar transform between minimal surfaces in $S^3$. In \cite{DVV}, the first named author of the present paper, Van der Veken and Vrancken define two transforms, called $(\pm)$transforms, between non-conformal harmonic surfaces in $S^3$, which is an extension of the polar transform. 

In \cite{S2}, the second named author of the present paper gives ($\pm$)transforms from a non-conformal harmonic surface $f$ in $R^3$, to holomorphic maps $f^{\pm}$ into $S^2$. It is as an extension of the classical fact that the Gauss map of a minimal surface in $R^3$ is holomorphic. \\

{\bf Remark 1.1.} In fact, in \cite{S2}, we showed that $f^{\pm}$ are holomorphic or anti-holomorphic. But in this paper, we will see that they are simultaneously holomorphic. \\

Let us consider the converse problem to \cite{S2}. Namely, for two distinct holomorphic maps $g, h: M\rightarrow S^2$ where $M$ is a Riemann surface, does there exist a non-conformal harmonic surface $f: M\rightarrow R^3$ such that $f^{+} = g$ and $f^{-} = h$? We solve this problem, and as a result, we get a representation formula for non-conformal harmonic surfaces in $R^3$, which is different from and simpler than that of Alarcon and Lopez \cite{AL}. 

\begin{thm}
Let $M$ be a simply connected Riemann surface with local complex coordinate $z$, and let $g, h: M\rightarrow S^2$ be holomorphic maps from $M$ to $S^2$, so that $g-h\neq 0$, $g+h\neq 0$. Then 
\[f = 4\ \mbox{Re}\int \frac{1}{|g-h|^2}\{(g\times h)-i(g-h)\}dz\]
gives a non-conformal harmonic surface in $R^3$ such that $f^{+} = g$ and $f^{-} = h$. 

Conversely, any non-conformal harmonic surface in $R^3$ can be locally represented as above. 
\end{thm}

{\bf Remark 1.2.} See \cite{BV}, \cite{AV} and \cite{S1} for $(\pm)$transforms in higher codimensional cases.

\section{Preliminaries}

Let $f:M\rightarrow R^3$ be a harmonic immersion of a Riemann surface $M$ into $R^3$ which is non-conformal at any point. We choose a local complex coordinate $z = x+i y$ on $M$. Then $f_{z\bar{z}} = 0$ by the harmonicity. We use the notation $\langle\ ,\ \rangle$ as the complex bilinear extension of the Euclidean metric. Since $\langle f_z, f_z\rangle$ is a non-zero holomorphic function, there exists a local complex coordinate $z$ such that $\langle f_z, f_z\rangle = -1$, which is called the adapted coordinate. It is equivalent to that
\[|f_y|^2-|f_x|^2 = 4,\ \ \ \ \langle f_x, f_y\rangle = 0.\]
Then there exists a smooth positive function $\phi$ such that
\[|f_x| = 2\sinh\phi,\ \ \ \ |f_y| = 2\cosh\phi.\]
Let $N$ be a unit normal vector field. 

Let $\theta$ be a smooth function so that
\[\cos\theta = \frac{|f_x|}{|f_y|} = \tanh\phi,\ \ \ \ \sin\theta = \mbox{sech}\ \phi.\]
Slanting the Gauss map, we define the $(\pm)$transforms $f^{\pm}: M\rightarrow S^2$ as follows: 
\[f^{\varepsilon} = \varepsilon(\sin\theta)\frac{f_y}{|f_y|}+(\cos\theta)N,\]
where $\varepsilon = +$ or $-$. In \cite{S2}, we showed that $f^{\pm}: M\rightarrow S^2$ are holomorphic or anti-holomorphic.

\section{Proof}

{\it Proof of Theorem 1.1.} (i) Let $f:M\rightarrow R^3$ be a non-conformal harmonic surface in $R^3$ with adapted coordinate. Using the notations in Section 2, we have
\[f^{+}+f^{-} = 2(\cos\theta)N,\ \ \ \ f^{+}-f^{-} = 2(\sin\theta)\frac{f_y}{|f_y|},\]
\[|f^{+}+f^{-}| = 2\cos\theta,\ \ \ \ |f^{+}-f^{-}| = 2\sin\theta,\]
and
\[N = \frac{f^{+}+f^{-}}{|f^{+}+f^{-}|},\ \ \ \ \frac{f_y}{|f_y|} = \frac{f^{+}-f^{-}}{|f^{+}-f^{-}|}.\]
From
\begin{eqnarray}
|f_y| = 2\cosh\phi = \frac{2}{\sin\theta} = \frac{4}{|f^{+}-f^{-}|}, \label{1} 
\end{eqnarray}
we have
\begin{eqnarray}
f_y = \frac{4}{|f^{+}-f^{-}|^2}(f^{+}-f^{-}). \label{2} 
\end{eqnarray}
Since
\begin{eqnarray}
|f_x| = |f_y|\cos\theta = \frac{2\cos\theta}{\sin\theta} = \frac{2|f^{+}+f^{-}|}{|f^{+}-f^{-}|}, \label{3} 
\end{eqnarray}
and
\[\frac{f_x}{|f_x|} = \frac{f_y}{|f_y|}\times N = \frac{2}{|f^{+}+f^{-}| |f^{+}-f^{-}|}(f^{+}\times f^{-}),\]
we get
\begin{eqnarray}
f_x = \frac{4}{|f^{+}-f^{-}|^2}(f^{+}\times f^{-}). \label{4} 
\end{eqnarray}
Hence
\[f_z = \frac{2}{|f^{+}-f^{-}|^2}\{(f^{+}\times f^{-})-i(f^{+}-f^{-})\},\]
which is holomorphic, and
\[f = 4\ \mbox{Re}\int \frac{1}{|f^{+}-f^{-}|^2}\{(f^{+}\times f^{-})-i(f^{+}-f^{-})\}dz.\]

Now we will show that $f^{+}$ and $f^{-}$ are simultaneously holomorphic. We denote the stereographic projection by $\Pi: S^2-\{(0,0,1)\}\rightarrow R^2$. Locally, we may write $f^{+} = \Pi^{-1}(P(z))$ and $f^{-} = \Pi^{-1}(Q(z))$ for some functions $P(z)$ and $Q(z)$. Let $P_1 = \mbox{Re}(P)$, $P_2 = \mbox{Im}(P)$, $Q_1 = \mbox{Re}(Q)$ and $Q_2 = \mbox{Im}(Q)$. Then we have 
\[f^{+} = \frac{1}{|P|^2+1}(2P_1, 2P_2, |P|^2-1),\ \ \ \ f^{-} = \frac{1}{|Q|^2+1}(2Q_1, 2Q_2, |Q|^2-1).\]
Since
\[f^{+}-f^{-} = \frac{2}{(|P|^2+1)(|Q|^2+1)} \left(\begin{array}{c}
P_1(|Q|^2+1)-Q_1(|P|^2+1) \\
P_2(|Q|^2+1)-Q_2(|P|^2+1) \\
|P|^2-|Q|^2
\end{array}\right),\]
and
\[|f^{+}-f^{-}|^2 = \frac{4|P-Q|^2}{(|P|^2+1)(|Q|^2+1)}, \]
we have by (\ref{2}), 
\[f_y = \frac{2}{|P-Q|^2} \left(\begin{array}{c}
P_1(|Q|^2+1)-Q_1(|P|^2+1) \\
P_2(|Q|^2+1)-Q_2(|P|^2+1) \\
|P|^2-|Q|^2
\end{array}\right).\]
From
\[f^{+}\times f^{-} = \frac{2}{(|P|^2+1)(|Q|^2+1)}
\left(\begin{array}{c}
P_2(|Q|^2-1)-Q_2(|P|^2-1) \\
Q_1(|P|^2-1)-P_1(|Q|^2-1) \\
2(P_1 Q_2-P_2 Q_1)
\end{array}\right),\]
we have by (\ref{4}), 
\[f_x = \frac{2}{|P-Q|^2} \left(\begin{array}{c}
P_2(|Q|^2-1)-Q_2(|P|^2-1) \\
Q_1(|P|^2-1)-P_1(|Q|^2-1) \\
2(P_1 Q_2-P_2 Q_1)
\end{array}\right).\]
Then we can get
\[f_z = \frac{1}{|P-Q|^2}\left(i (\bar{P}-\bar{Q})(PQ-1), (\bar{P}-\bar{Q})(PQ+1), -i (\bar{P}-\bar{Q})(P+Q) \right) \]
\[= \left(i\frac{PQ-1}{P-Q}, \frac{PQ+1}{P-Q}, -i\frac{P+Q}{P-Q} \right).\]
Since $f_z$ is holomorphic, both $P(z)$ and $Q(z)$ are holomorphic. Hence, both $f^{+}$ and $f^{-}$ are holomorphic, as desired. 

(ii) Let $M$ be a simply connected Riemann surface with local complex coordinate $z = x+i y$. Let $g, h:M\rightarrow S^2$ be holomorphic maps from $M$ to $S^2$, so that $g-h \neq 0$ and $g+h \neq 0$. We will show that 
\[\frac{1}{|g-h|^2}\{(g\times h)-i(g-h)\} :M \rightarrow C^3 \]
is holomorphic. Using the stereographic projection $\Pi$, we may write locally $g = \Pi^{-1}(P(z))$ and $h = \Pi^{-1}(Q(z))$ for some holomorphic functions $P(z)$ and $Q(z)$. Let $P_1 = \mbox{Re}(P)$, $P_2 = \mbox{Im}(P)$, $Q_1 = \mbox{Re}(Q)$ and $Q_2 = \mbox{Im}(Q)$. Then we have 
\[g = \frac{1}{|P|^2+1}(2P_1, 2P_2, |P|^2-1),\ \ \ \ h = \frac{1}{|Q|^2+1}(2Q_1, 2Q_2, |Q|^2-1).\]
By the same computation as in the part (i), we can see that
\[\frac{1}{|g-h|^2}\{(g\times h)-i(g-h)\} = \frac{1}{2}\left(i\frac{PQ-1}{P-Q}, \frac{PQ+1}{P-Q}, -i\frac{P+Q}{P-Q} \right),\]
which is holomorphic. 

Then
\[f = 4\ \mbox{Re}\int \frac{1}{|g-h|^2}\{(g\times h)-i(g-h)\}dz \]
gives a harmonic map into $R^3$. We have
\[f_x = \frac{4}{|g-h|^2}(g\times h),\ \ \ \ f_y = \frac{4}{|g-h|^2}(g-h), \]
\[\langle f_x, f_y\rangle = 0, \ \ \ \ |f_y|^2 = \frac{16}{|g-h|^2}. \]
Noting that
\[|g\times h|^2 = 1-\langle g, h\rangle^2 = (1+\langle g, h\rangle)(1-\langle g, h\rangle) = \frac{1}{4}|g+h|^2 |g-h|^2 ,\]
we have
\[|f_x|^2 = \frac{4|g+h|^2}{|g-h|^2}, \]
and $|f_y|^2-|f_x|^2 = 4$. So $f$ gives a non-conformal harmonic surface in $R^3$ with adapted coordinate. 

We can compute that 
\[f_x\times f_y = \frac{16}{|g-h|^4}\{(g\times h)\times(g-h)\} = \frac{8}{|g-h|^2}(g+h),\]
and
\[N = \frac{g+h}{|g+h|}.\]
Let $\theta$ be a smooth function so that
\[\cos\theta = \frac{|f_x|}{|f_y|} = \frac{1}{2}|g+h|,\ \ \ \ \sin\theta = \frac{1}{2}|g-h|,\]
where we note that $|g+h|^2+|g-h|^2 = 4$. Then we can get $f^{+} = g$ and $f^{-} = h$. 

By (i) and (ii), we have proved the theorem. 
\\

From the proof of Theorem 1.1, we also get the following: 

\begin{cor}
Any non-conformal harmonic surface in $R^3$ can be locally represented by
\[f = 2\ \mbox{Re}\int \left(i\frac{PQ-1}{P-Q}, \frac{PQ+1}{P-Q}, -i\frac{P+Q}{P-Q} \right) dz, \]
where $P(z)$ and $Q(z)$ are holomorphic functions. Here we have $f^{+} = \Pi^{-1}\circ P$ and $f^{-} = \Pi^{-1}\circ Q$. 
\end{cor}

{\bf Remark 3.1.} From the proof of Theorem 1.1, we can see that 
\[|g-h|^2 = \frac{4|P-Q|^2}{(|P|^2+1)(|Q|^2+1)}, \ \ \ \ |g+h|^2 = \frac{4 |\bar{P}Q+1|^2}{(|P|^2+1)(|Q|^2+1)}. \]
So $g-h \neq 0$ and $g+h \neq 0$ are equivalent to that $P \neq Q$ and $\bar{P}Q \neq -1$. \\

Let $f$ be a non-conformal harmonic surface in $R^3$ with constant $f^{+}$ or $f^{-}$. By identification, we consider the case of constant $f^{-}$. Up to congruence, we may assume that $f^{-} = (0, 0, -1)$ and $Q = 0$ in Corollary 3.1. Then $f$ can be locally represented by 
\[f = 2\ \mbox{Re}\int\left(\frac{-i}{P}, \frac{1}{P}, -i \right)dz\]
for a holomorphic function $P(z)$. Set
\[\psi = \int\frac{2}{P}dz,\]
which is holomorphic, and $P = 2/\psi'$. Then we have
\[f = (\mbox{Im}(\psi), \mbox{Re}(\psi), 2y). \]

\begin{cor}
Up to congruence, any non-conformal harmonic surface $f$ in $R^3$ with constant $f^{-}$ can be locally represented by 
\[f = (\mbox{Im}(\psi), \mbox{Re}(\psi), 2y), \]
where $\psi(z)$ is a holomorphic function with $\psi'(z) \neq 0$. 
\end{cor}

\section{Completeness and curvature}

For a non-conformal harmonic surface $f$ as in Corollary 3.1, by (\ref{1}), (\ref{3}) and Remark 3.1, we have
\[|f_x|^2 = \frac{4|\bar{P}Q+1|^2}{|P-Q|^2},\]
\[|f_y|^2 = \frac{4(|P|^2+1)(|Q|^2+1)}{|P-Q|^2},\]
and the first fundamental form 
\[ds^2 = 4\left(\frac{|\bar{P}Q+1|^2}{|P-Q|^2}dx^2+\frac{(|P|^2+1)(|Q|^2+1)}{|P-Q|^2}dy^2\right). \]
Noting that
\[\frac{(|P|^2+1)(|Q|^2+1)}{|P-Q|^2} = \frac{|\bar{P}Q+1|^2}{|P-Q|^2}+1,\]
we get the following completeness criterion. 

\begin{thm}
Let $f:C\rightarrow R^3$ be a non-conformal harmonic surface represented as in Corollary 3.1, where the domain is the complex plane $C$. Suppose that $P(z), Q(z)$ are holomorphic functions on $C$, $P \neq Q$, $\bar{P}Q \neq -1$, and there exists a positive constant $\delta$ such that
\[\frac{|\bar{P}Q+1|^2}{|P-Q|^2} \geq \delta. \]
Then the surface $f$ is complete. 
\end{thm}

For a non-conformal harmonic surface $f$ as in Corollary 3.1, the Gaussian curvature $K$ is given by
\[K = -\frac{4|\langle f_{zz}, N \rangle|^2}{|f_x|^2|f_y|^2}. \]
We have
\[f_z = \left(i\frac{PQ-1}{P-Q}, \frac{PQ+1}{P-Q}, -i\frac{P+Q}{P-Q}\right), \]
\[N = \frac{g+h}{|g+h|} = \hspace{8cm} \]
\[\frac{1}{2|\bar{P}Q+1|\sqrt{(|P|^2+1)(|Q|^2+1)}}
\left( \begin{array} {c} 
(P+\bar{P})(|Q|^2+1)+(Q+\bar{Q})(|P|^2+1) \\
-i \{(P-\bar{P})(|Q|^2+1)+(Q-\bar{Q})(|P|^2+1)\} \\
2(|P|^2|Q|^2-1) 
\end{array} \right), \]
\[f_{zz} = \frac{1}{(P-Q)^2}
\left( \begin{array} {c} 
i \{Q'(P^2-1)-P'(Q^2-1)\} \\
Q'(P^2+1)-P'(Q^2+1) \\
-2i(PQ'-P'Q)
\end{array} \right), \]
and
\[\langle f_{zz}, N \rangle = i\frac{(|Q|^2+1)(\bar{P}Q+1)P'+(|P|^2+1)(P\bar{Q}+1)Q'}{(P-Q)|\bar{P}Q+1|\sqrt{(|P|^2+1)(|Q|^2+1)}}. \]
Thus we get the following:

\begin{prop}
Let $f:M\rightarrow R^3$ be a non-conformal harmonic surface represented as in Corollary 3.1. Then the Gaussian curvature $K$ is given by
\[K = -\frac{|P-Q|^2 \left|(|Q|^2+1)(\bar{P}Q+1)P'+(|P|^2+1)(P\bar{Q}+1)Q' \right|^2}{4|\bar{P}Q+1|^4 (|P|^2+1)^2 (|Q|^2+1)^2}. \]
\end{prop}

\vspace{2mm}

From this proposition, we find that
\[|K| \leq \frac{|P-Q|^2}{2|\bar{P}Q+1|^2}\left(\frac{|P'|^2}{(|P|^2+1)^2}+\frac{|Q'|^2}{(|Q|^2+1)^2}\right), \]
and noting that
\[|dg|^2 = |g_x|^2+|g_y|^2 = \frac{8|P'|^2}{(|P|^2+1)^2},\ \ \ \ |dh|^2 = |h_x|^2+|h_y|^2 = \frac{8|Q'|^2}{(|Q|^2+1)^2}, \]
we have
\[|K| \leq \frac{|g-h|^2}{16|g+h|^2}(|dg|^2+|dh|^2). \]
Let $dV$ be the area element of $M$ induced by $f$, which is given by
\[dV = \frac{8|g+h|}{|g-h|^2}dxdy.\]
So if 
\[\int_{M}\frac{|dg|^2+|dh|^2}{|g+h|}dxdy < \infty, \]
then $f$ has finite total curvature.

\section{Examples}

{\bf Example 5.1.} In Corollary 3.2, set $\psi(z) = z^2$. Then we have
\[f(x, y) = (2xy, x^2-y^2, 2y), \ \ \ \ (x, y) \neq (0,0). \]
It is not complete, and $f^{-} = (0,0,-1)$. \\

{\bf Example 5.2.} In Corollary 3.2, set $\psi(z) = iae^{-iz}$, where $a$ is a positive real number. Then
\[f(x, y) = (ae^y\cos{x}, ae^y\sin{x}, 2y).\]
It is a complete rotational surface, and $f^{-} = (0,0,-1)$ (cf. \cite{S2}). \\

{\bf Example 5.3.} In Corollary 3.1, set $P(z) = ae^z$ and $Q(z) = e^z/a$, where $a$ is a real number such that $a \neq 0, \pm 1$. Then we can see that 
\[\frac{|\bar{P}Q+1|^2}{|P-Q|^2} = \frac{a^2}{(a^2-1)^2}\cdot\frac{(e^{2x}+1)^2}{e^{2x}} \geq \frac{4a^2}{(a^2-1)^2}. \]
By Theorem 4.1, the surface $f$ is complete. In fact we have
\[f(x, y) = \frac{2}{a^2-1} \left( \begin{array}{c}
-2a\sinh{x} \sin{y} \\
2a\sinh{x} \cos{y} \\
(a^2+1)y
\end{array} \right).\]
It is a helicoid, and gives a minimal surface with non-conformal harmonic parametrization. \\

{\bf Example 5.4.} In Corollary 3.1, set $P(z) = ae^{iz}$ and $Q(z) = e^{iz}/a$, where $a$ is a real number such that $a \neq 0, \pm 1$. Then 
\[\frac{|\bar{P}Q+1|^2}{|P-Q|^2} = \frac{a^2}{(a^2-1)^2}\cdot\frac{(1+e^{2y})^2}{e^{2y}} \geq \frac{4a^2}{(a^2-1)^2}. \]
By Theorem 4.1, the surface $f$ is complete. In fact we have
\[f(x, y) = \frac{2}{a^2-1} \left( \begin{array}{c}
2a\cosh{y} \cos{x} \\
2a\cosh{y} \sin{x} \\
(a^2+1)y
\end{array} \right). \]
It is a catenoid-like surface, but does not give a minimal surface. \\

{\bf Remark 5.1.} Let $f:C\rightarrow R^3$ be a non-conformal harmonic surface, where the domain is the complex plane $C$. By the small Picard theorem, if $f^{+}$ (resp. $f^{-}$) omits three points on $S^2$, then $f^{+}$ (resp. $f^{-}$) has to be constant. A Liouville type theorem was stated in \cite{S2}. Examples 5.3 and 5.4 give non-conformal harmonic surfaces $f:C\rightarrow R^3$ such that both $f^{+}$ and $f^{-}$ omit just two points on $S^2$.

\section{Remark on the distortions}

Let $f:M\rightarrow R^3$ be a harmonic immersion of a Riemann surface $M$ into $R^3$. Let $z = x+iy$ be a local complex coordinate on $M$, and $N$ the unit normal vector field. The distortion function $D_f$ of $f$ is given by
\[D_f = \frac{|f_x|^2+|f_y|^2}{2|f_x\times f_y|}. \]
The distortion function $D_N$ of the Gauss map is given by
\[D_N = \frac{|N_x|^2+|N_y|^2}{2|N_x\times N_y|}, \]
when the Gauss map is regular. In \cite{K} Kalaj showed that $D_N = D_f$, provided that the Gauss map is regular. But the proof is a little complicated. So, in the following, we give a simple proof. 

\vspace{2mm}

Let
\[E = |f_x|^2,\ \ \ \ F = \langle f_x, f_y\rangle,\ \ \ \ G = |f_y|^2,\]
\[\ell = \langle f_{xx}, N\rangle,\ \ \ \ m = \langle f_{xy}, N\rangle,\ \ \ \ -\ell = \langle f_{yy}, N\rangle, \]
which are the components of the first and second fundamental forms.

(i) Around a non-conformal point, we choose an adapted coordinate as in Section 2. Then $F = 0, G = E+4$, and
\[D_f = \frac{E+G}{2\sqrt{EG}}.\]
We have
\[N_x = -\frac{\ell}{E}f_x-\frac{m}{G}f_y,\ \ \ \ N_y = -\frac{m}{E}f_x+\frac{\ell}{G}f_y, \]
and
\[|N_x|^2 = \frac{\ell^2}{E}+\frac{m^2}{G},\ \ \ \ |N_y|^2 = \frac{m^2}{E}+\frac{\ell^2}{G},\]
\[|N_x\times N_y| = \frac{\ell^2+m^2}{\sqrt{EG}}.\]
Noting that $\ell^2+m^2 > 0$ by the regularity of the Gauss map, we get
\[D_N = \frac{E+G}{2\sqrt{EG}}.\]

(ii) At a conformal point $p$, we have $E = G, F = 0$, and $D_f = 1$. At $p$, we have
\[N_x = -\frac{\ell}{E}f_x-\frac{m}{E}f_y,\ \ \ \ N_y = -\frac{m}{E}f_x+\frac{\ell}{E}f_y, \]
and
\[|N_x|^2 = |N_y|^2 = |N_x\times N_y| = \frac{\ell^2+m^2}{E}. \]
Since $\ell^2+m^2 > 0$ by the regularity of the Gauss map, we get $D_N = 1$ at $p$. 

By (i) and (ii), we have shown that $D_N = D_f$, provided that the Gauss map is regular.

\vspace{5mm}

Bart Dioos, KU Leuven, Tutorial services, Faculty of Engineering Science, Celestijnenlaan 200i - box 2201, 3001 Leuven, Belgium 

E-mail: bart.dioos@kuleuven.be 

\vspace{3mm}

Makoto Sakaki, Graduate School of Science and Technology, Hirosaki University, Hirosaki 036-8561, Japan 

E-mail: sakaki@hirosaki-u.ac.jp 

\end{document}